\documentclass{amsart}
\usepackage{amssymb}
\newtheorem{thm}{Theorem}
\newtheorem{lem}{Lemma}
\newtheorem{prop}{Proposition}

\newtheorem{df}{Definition}

\begin{document}

\bibliographystyle{plain}

\title[Milo\v s Arsenovi\'c and Romi F. Shamoyan]{Sharp theorems on multipliers and distances in harmonic function spaces in higher dimension}

\author[]{Milo\v s Arsenovi\' c$\dagger$}
\author[]{Romi F. Shamoyan}

\address{Faculty of mathematics, University of Belgrade, Studentski Trg 16, 11000 Belgrade, Serbia}
\email{\rm arsenovic@matf.bg.ac.rs}

\address{Bryansk University, Bryansk Russia}
\email{\rm rshamoyan@yahoo.com}

\thanks{$\dagger$ Supported by Ministry of Science, Serbia, project M144010}

\date{}

\begin{abstract}
We present new sharp results concerning multipliers and distance estimates in various spaces of harmonic functions in the unit ball of $\mathbb R^n$.
\end{abstract}

\maketitle

\footnotetext[1]{Mathematics Subject Classification 2010 Primary 42B15, Secondary 42B30.  Key words
and Phrases: Multipliers, harmonic functions, Bergman spaces, mixed norm spaces, distance estimates.}

\section{Introduction and preliminaries}

The aim of this paper is twofold. One is to describe spaces of multipliers between certain spaces of harmonic functions on the unit ball. We note that so far there are no results in this direction in the multidimensional case, where the use of spherical harmonics is a natural substitute for power series expansion. In fact, even the case of the unit disc has not been extensively studied in this context. We refer the reader to \cite{SW}, where multipliers between harmonic Bergman type classes were considered, and to \cite{Pa1} and \cite{Pa2} for the case of harmonic Hardy classes. Most of
our results are present in these papers in the special case of the unit disc.

The other topic we investigate is distance estimates in spaces of harmonic functions on the unit ball. This line
of investigation can be considered as a continuation of papers \cite{AS1}, \cite{SM1} and \cite{SM2}.

Let $\mathbb B$ be the open unit ball in $\mathbb R^n$, $\mathbb S = \partial \mathbb B$ is the unit sphere in $\mathbb R^n$, for $x \in \mathbb R^n$ we have $x = rx'$, where $r = |x| = \sqrt{\sum_{j=1}^n x_j^2}$ and $x' \in \mathbb S$. Normalized Lebesgue measure on $\mathbb B$ is denoted by $dx = dx_1 \ldots dx_n = r^{n-1}dr dx'$ so that $\int_{\mathbb B} dx = 1$. We denote the space of all harmonic functions in an open set $\Omega$ by $h(\Omega)$.

We consider harmonic weighted Bergman spaces $A^p_\alpha (\mathbb B)$ on $\mathbb B$ defined for $\alpha > -1$ and
$0 < p \leq \infty$ by
$$A^p_\alpha(\mathbb B) = \left\{ f \in h(\mathbb B) : \| f \|_{p, \alpha} = \left( \int_{\mathbb B}
|f(rx')|^p(1-r)^\alpha r^{n-1}dr dx' \right)^{1/p} < \infty \right\},$$
$$A^\infty_\alpha(\mathbb B) = \left\{ f \in h(\mathbb B) : \|f\|_{\infty, \alpha} = \sup_{x \in \mathbb B}
|f(x)|(1-|x|)^\alpha < \infty \right\}.$$
It is easy to show that the spaces $A^p_\alpha = A^p_\alpha (\mathbb B)$ are Banach spaces for $1 \leq p \leq \infty$ and complete metric spaces for $0<p<1$.

For $0<p<\infty$, $0 \leq r < 1$ and $f \in h(\mathbb B)$ we set
$$M_p(f, r) = \left( \int_{\mathbb S} |f(rx')|^p dx' \right)^{1/p},$$
with the usual modification to cover the case $p = \infty$. Weighted Hardy spaces are defined, for $\alpha \geq 0$ and
$0 < p \leq \infty$, by
$$H^p_\alpha (\mathbb B) = H^p_\alpha = \{ f \in h(\mathbb B) : \| f \|_{p, \alpha} = \sup_{r<1} M_p(f, r)(1-r)^\alpha
< \infty \}.$$
For $\alpha = 0$ the space $H^p_\alpha$ is denoted simply by $H^p$.

For $0 < p \leq \infty$, $0 < q \leq \infty$ and $\alpha > 0$ and we consider mixed (quasi)-norms
$\| f \|_{p,q;\alpha}$ defined by
\begin{equation}\label{qpnorm}
\| f \|_{p,q;\alpha} = \left( \int_0^1 M_q(f, r)^p (1-r^2)^{\alpha p - 1} r^{n-1}dr \right)^{1/p}, \qquad
f \in h(\mathbb B),
\end{equation}
again with the usual interpretation for $p = \infty$, and the corresponding spaces
$$B^{p,q}_\alpha(\mathbb B) = B^{p,q}_\alpha = \{ f \in h(\mathbb B) : \| f \|_{p,q;\alpha} < \infty \}.$$
It is not hard to show that these spaces are complete metric spaces and that for $\min(p,q) \geq 1$ they are Banach spaces.

Note that $A^\infty_\alpha = H^\infty_\alpha$ for $\alpha \geq 0$ and $B^{\infty, q}_\alpha = H^q_\alpha$ for $0 < q \leq \infty$, $\alpha > 0$. We also have, for $0 < p_0 \leq p_1 \leq \infty$, $B^{p_0, 1}_\alpha \subset B^{p_1, 1}_\alpha$, see \cite{DS}.

Next we need certain facts on spherical harmonics and Poisson kernel, see \cite{SW} for a detailed exposition. Let $Y^{(k)}_j$ be the spherical harmonics of order $k$, $j \leq 1 \leq d_k$, on $\mathbb S$. Next,
$$Z_{x'}^{(k)}(y') = \sum_{j=1}^{d_k} Y_j^{(k)}(x') \overline{Y_j^{(k)}(y')}$$
are zonal harmonics of order $k$. Note that the spherical harmonics $Y^{(k)}_j$, ($k \geq 0$, $1 \leq j \leq d_k$) form an orthonormal basis of $L^2(\mathbb S, dx')$. Every $f \in h(\mathbb B)$ has an expansion
$$f(x) = f(rx') = \sum_{k=0}^\infty r^k b_k\cdot Y^k(x'),$$
where $b_k = (b_k^1, \ldots, b_k^{d_k})$, $Y^k = (Y_1^{(k)}, \ldots, Y_{d_k}^{(k)})$ and $b_k\cdot Y^k$ is interpreted in the scalar product sense: $b_k \cdot Y^k = \sum_{j=1}^{d_k} b_k^j Y_j^{(k)}$. We often write, to stress dependence on a function $f \in h(\mathbb B)$, $b_k = b_k(f)$ and $b_k^j = b_k^j(f)$, in fact we have linear functionals $b_k^j$,
$k \geq 0, 1 \leq j \leq d_k$ on the space $h(\mathbb B)$.

We denote the Poisson kernel for the unit ball by $P(x, y')$, it is given by
\begin{align*}
P(x, y') = P_{y'}(x) & = \sum_{k=0}^\infty r^k \sum_{j=1}^{d_k} Y^{(k)}_j(y') Y^{(k)}_j(x') \\
& = \frac{1}{n\omega_n} \frac{1-|x|^2}{|x-y'|^n}, \qquad x = rx' \in \mathbb B, \quad y' \in \mathbb S,
\end{align*}
where $\omega_n$ is the volume of the unit ball in $\mathbb R^n$. We are going to use also a Bergman kernel for $A^p_\beta$ spaces, this is the following function
\begin{equation}\label{bker}
Q_\beta(x, y) = 2 \sum_{k=0}^\infty \frac{\Gamma(\beta + 1 + k + n/2)}{\Gamma(\beta + 1) \Gamma(k + n/2)}
r^k \rho^k Z_{x'}^{(k)}(y'), \qquad x = rx', \; y = \rho y' \in \mathbb B.
\end{equation}
For details on this kernel we refer to \cite{DS}, where the following theorem can be found.
\begin{thm}[\cite{DS}]\label{intrep}
Let $p \geq 1$ and $\beta \geq 0$. Then for every $f \in A^p_\beta$ and $x \in \mathbb B$ we have
$$f(x) = \int_0^1 \int_{\mathbb S^{n-1}} Q_\beta(x, y) f(\rho y')(1-\rho^2)^\beta \rho^{n-1}d\rho dy',
\qquad y = \rho y'.$$
\end{thm}
This theorem is a cornerstone for our approach to distance problems in the case of the unit ball. The following lemma from \cite{DS} gives estimates for this kernel.
\begin{lem}\label{DSlemma}
1. Let $\beta > 0$. Then, for $x = rx', y = \rho y' \in \mathbb B$ we have
$$|Q_\beta(x, y)| \leq C\frac{(1-r\rho)^{-\{ \beta \}}}{|\rho x - y'|^{n + [\beta]}} + \frac{C}{(1-r\rho)^{1+\beta}}.$$
If, moreover, $\beta > 0$ is an integer, then we have
$$|Q_\beta(x, y)| \leq \frac{C}{|\rho x - y'|^{n+\beta}}.$$
2. Let $\beta > -1$. Then
$$\int_{\mathbb S^{n-1}} |Q_\beta(rx', y)| dx' \leq \frac{C}{(1-r\rho)^{1+\beta}},
\qquad |y| = \rho, \quad 0 \leq r < 1.$$
3. Let $\beta > n-1$, , $0 \leq r < 1$ and $y'\in \mathbb S^{n-1}$. Then
$$\int_{\mathbb S^{n-1}} \frac{dx'}{|rx' - y'|^\beta} \leq \frac{C}{(1-r)^{\beta-n+1}}.$$
\end{lem}

\begin{lem}[\cite{DS}]\label{rro}
Let $\alpha > -1$ and $\lambda > \alpha + 1$. Then
$$\int_0^1 \frac{(1-r)^\alpha}{(1-r\rho)^\lambda} dr \leq C (1-\rho)^{\alpha + 1 - \lambda}, \qquad 0 \leq \rho < 1.$$
\end{lem}

\begin{lem}\label{wellkn}
Let $G(r)$, $0 \leq r < 1$, be a positive increasing function. Then, for $\alpha > -1$, $\beta > -1$, $\gamma \geq 0$ and $0 < q \leq 1$ we have
\begin{equation}\label{welleq}
\left( \int_0^1 G(r) \frac{(1-r)^\beta}{(1- \rho r)^\gamma} r^\alpha dr \right)^q \leq C
\int_0^1 G(r)^q \frac{(1-r)^{\beta q + q - 1}}{(1-\rho r)^{q\gamma}} r^\alpha dr, \quad 0 \leq \rho < 1.
\end{equation}
\end{lem}

A special case of the above lemma appears in \cite{AS}, for reader's convenience we produce a proof.

{\it Proof.} We use a subdivision of $I = [0, 1)$ into subintervals $I_k = [r_k, r_{k+1})$, $k \geq 0$, where
$r_k = 1-2^{-k}$. Since $1- \rho r_k \asymp 1- \rho r_{k+1}$, $0 \leq \rho < 1$, we have
\begin{align*}
J & = \left( \int_0^1 G(r) \frac{(1-r)^\beta}{(1- \rho r)^\gamma} r^\alpha dr \right)^q = \left( \sum_{k\geq 0}
\int_{I_k}  G(r) \frac{(1-r)^\beta}{(1- \rho r)^\gamma} r^\alpha dr \right)^q \\
& \leq \sum_{k\geq 0} \left( \int_{I_k}  G(r) \frac{(1-r)^\beta}{(1- \rho r)^\gamma} r^\alpha dr \right)^q \leq
\sum_{k \geq 0} 2^{-kq\beta} G^q(r_{k+1}) \left( \int_{I_k} \frac{r^\alpha dr}{(1-\rho r)^\gamma} \right)^q \\
& \leq C \sum_{k\geq 0} 2^{-kq\beta} G^q(r_{k+1}) 2^{-kq} (1- \rho r_{k+1})^{-q\gamma}\\
& \leq C \sum_{k\geq 0} 2^{-kq\beta} G^q(r_{k+1}) 2^{-kq}(1- \rho r_k)^{-q\gamma} \\
& \leq C \sum_{k\geq 0} G^q(r_{k+1}) \int_{I_{k+1}} \frac{(1-r)^{\beta q + q - 1} r^\alpha dr}{(1-\rho r)^{q\gamma}} \\
& \leq C \int_0^1 G(r)^q \frac{(1-r)^{\beta q + q - 1}}{(1-\rho r)^{q\gamma}} r^\alpha dr. \qquad \Box
\end{align*}

\begin{lem}\label{qbeta}
For $\delta > -1$, $\gamma > n + \delta$ and integer $\beta > 0$ we have
$$\int_{\mathbb B} |Q_\beta(x, y)|^{\frac{\gamma}{n+\beta}} (1-|y|)^\delta dy \leq C (1-|x|)^{\delta - \gamma + n},
\qquad x \in \mathbb B.$$
\end{lem}

{\it Proof.} Using Lemma \ref{DSlemma} and Lemma \ref{rro} we obtain:
\begin{align*}
\int_{\mathbb B} |Q_\beta(x, y)|^{\frac{\gamma}{n+\beta}} (1-|y|)^\delta dy & \leq C \int_{\mathbb B}
\frac{(1-|y|)^\delta}{|\rho rx' - y'|^\gamma} dy \\
& \leq C \int_0^1 (1-\rho)^\delta \int_{\mathbb S} \frac{dy'}{|\rho r x' - y' |^\gamma} dy' d\rho \\
& \leq C \int_0^1 (1-\rho)^\delta (1-r\rho)^{n - \gamma - 1}d\rho \\
& \leq C (1-r)^{n + \delta - \gamma}. \qquad \Box
\end{align*}

\begin{lem}[\cite{DS}]\label{gamma}
For real $s, t$ such that $s>-1$ and $2t+n>0$ we have
$$\int_0^1 (1-r^2)^s r^{2t+n-1} dr = \frac{1}{2} \frac{\Gamma(s+1)\Gamma(n/2+t)}{\Gamma(s+1+n/2+t)}.$$
\end{lem}

We set $\mathbb R^{n+1}_+ = \{ (x, t) : x \in \mathbb R^n, t > 0 \} \subset \mathbb R^{n+1}$. We usually denote the points in $\mathbb R^{n+1}_+$ by $z = (x, t)$ or $w = (y, s)$ where $x, y \in \mathbb R^n$ and $s, t > 0$. 

For $0 < p < \infty$ and $\alpha > -1$ we consider spaces
$$\tilde A^p_\alpha(\mathbb R^{n+1}_+) = \tilde A^p_\alpha = \left\{ f \in h(\mathbb R_+^{n+1}) : \int_{\mathbb R_+^{n+1}} |f(x, t)|^p t^\alpha dx dt < \infty \right\}.$$
Also, for $p = \infty$ and $\alpha > 0$, we set
$$\tilde A^\infty_\alpha(\mathbb R^{n+1}_+) = \tilde A^\infty_\alpha = \left\{ f \in h(\mathbb R_+^{n+1}) :
\sup_{(x, t) \in \mathbb R^{n+1}_+} |f(x, t)|t^\alpha < \infty \right\}.$$
These spaces have natural (quasi)-norms, for $1 \leq p \leq \infty$ they are Banach spaces and for $0<p\leq 1$ they are
complete metric spaces.

We denote the Poisson kernel for $\mathbb R^{n+1}_+$ by $P(x, t)$, i.e.
$$P(x, t) = c_n \frac{t}{(|x|^2 + t^2)^{\frac{n+1}{2}}}, \qquad x \in \mathbb R^n, t > 0.$$
For an integer $m \geq 0$ we introduce a Bergman kernel $Q_m(z, w)$, where $z = (x, t) \in \mathbb R^{n+1}_+$ and $w = (y, s) \in \mathbb R^{n+1}_+$, by
$$Q_m(z, w) = \frac{(-2)^{m+1}}{m!} \frac{\partial^{m+1}}{\partial t^{m+1}} P(x-y, t+s).$$
The terminology is justified by the following result from \cite{DS}.

\begin{thm}\label{brthm}
Let $0<p<\infty$ and $\alpha > -1$. If $0<p\leq 1$ and $m \geq \frac{\alpha +n +1}{p} - (n+1)$ or $1\leq p < \infty$
and $m> \frac{\alpha +1}{p} - 1$, then
\begin{equation}\label{brep}
f(z) = \int_{\mathbb R^{n+1}_+} f(w) Q_m(z, w) s^m dy ds, \qquad f \in \tilde A^p_\alpha,\quad z \in \mathbb R^{n+1}_+.
\end{equation}
\end{thm}

The following elementary estimate of this kernel is contained in \cite{DS}:
\begin{equation}\label{estq}
|Q_m(z, w)| \leq C \left[ |x-y|^2 + (s+t)^2 \right]^{-\frac{n+m+1}{2}}, \quad z = (x, t), w = (y, s) \in \mathbb R^{n+1}_+.
\end{equation}

\section{Multipliers on spaces of harmonic functions}

In this section we present our results on multipliers between spaces of harmonic functions on the unit ball. The
following definitions are needed to formulate these theorems.

\begin{df}
For a double indexed sequence of complex numbers
$$c = \{ c_k^j : k \geq 0, 1 \leq j \leq d_k \}$$
and a harmonic function $f(rx') = \sum_{k=0}^\infty \sum_{j=1}^{d_k} r^k b_k^j(f) Y^{(k)}_j(x')$ we define
$$(c \ast f) (rx') = \sum_{k=0}^\infty \sum_{j=1}^{d_k} r^k c_k^j b_k^j(f) Y^{(k)}_j(x'), \qquad rx' \in \mathbb B,$$
if the series converges in $\mathbb B$. Similarly we define convolution of $f, g \in h(\mathbb B)$ by
$$(f \ast g)(rx') = \sum_{k=0}^\infty \sum_{j=1}^{d_k} r^k b_k^j(f)b_k^j(f) Y_j^{(k)}(x'), \qquad rx'\in \mathbb B,$$
it is easily seen that $f \ast g$ is defined and harmonic in $\mathbb B$.
\end{df}

\begin{df}
For $t > 0$ and a harmonic function $f(x) = \sum_{k=0}^\infty b_k(f)Y^k(x')$ on the unit ball we define a
fractional derivative of order $t$ of $f$ by the following formula:
$$(\Lambda_t f)(x) = \sum_{k=0}^\infty r^k \frac{\Gamma(k+n/2 + t)}{\Gamma(k+n/2)\Gamma(t)}b_k(f)\cdot Y^k(x'),
\qquad x = rx' \in \mathbb B.$$
\end{df}
Clearly, for $f \in h(\mathbb B)$ and $t>0$ the function $\Lambda_t h$ is also harmonic in $\mathbb B$.

\begin{df}
Let $X$ and $Y$ be subspaces of $h(\mathbb B)$. We say that a double indexed sequence $c$ is a multiplier from $X$ to $Y$ if $c \ast f \in Y$ for every $f \in X$. The vector space of all multipliers from $X$ to $Y$ is denoted by
$M_H(X, Y)$.
\end{df}

Clearly every multiplier $c \in M_H(X, Y)$ induces a linear map $M_c : X \rightarrow Y$. If, in addition, $X$ and $Y$ are (quasi)-normed spaces such that all functionals $b_k^j$ are continuous on both spaces $X$ and $Y$, then the map
$M_c : X \rightarrow Y$ is continuous, as is easily seen using the Closed Graph Theorem. We note that this holds
for all spaces we consider in this paper : $A^p_\alpha$, $B^{p,q}_\alpha$ and $H^p_\alpha$.

The first part of the lemma below appeared, in dimension two, in \cite{Pa2}.

\begin{lem}\label{intcon}
Let $f, g \in h(\mathbb B)$ have expansions
\begin{equation*}
f(rx')  = \sum_{k=0}^\infty r^k \sum_{j=1}^{d_k} b_k^j Y^{(k)}_j(x'), \qquad
g(rx')  = \sum_{l=0}^\infty r^k \sum_{i=1}^{d_k} c_l^i Y^{(l)}_i(x').
\end{equation*}
Then we have
\begin{equation*}
\int_{\mathbb S} (g \ast P_{y'})(rx') f(\rho x') dx' = \sum_{k=0}^\infty r^k\rho^k \sum_{j=1}^{d_k} b_k^j c_k^j
Y^{(k)}_j(y'), \qquad y' \in \mathbb S, \quad 0 \leq r, \rho  < 1.
\end{equation*}
Moreover, for every $m > -1$, $y' \in \mathbb S$ and $0 \leq r, \rho < 1$ we have
\begin{equation*}
\int_{\mathbb S} (g \ast P_{y'})(rx') f(\rho x') dx' = 2 \int_0^1 \int_{\mathbb S} \Lambda_{m+1}(g \ast P_{y'})(rRx')
f(\rho Rx') (1-R^2)^m R^{n-1} dx' dR.
\end{equation*}
\end{lem}

{\it Proof.} The first assertion of this lemma easily follows from the orthogonality relations for spherical harmonics $Y^{(k)}_j$. Using Lemma \ref{gamma} and orthogonality relations we have
\begin{align*}
I & =  2 \int_0^1 \int_{\mathbb S} \Lambda_{m+1}(g \ast P_{y'})(rRx')
f(\rho Rx') (1-R^2)^m R^{n-1} dx' dR\\
& = 2\int_0^1 \sum_{k=0}^\infty r^k\rho^k R^{2k+n-1} (1-R^2)^m \frac{\Gamma(k+n/2 + m + 1)}{\Gamma(k+n/2)\Gamma(m+1)}\sum_{j=1}^{d_k} b_k^j c_k^j Y^{(k)}_j dR \\
& = \sum_{k=0}^\infty r^k \rho^k \sum_{j=1}^{d_k} b_k^j c_k^j Y^{(k)}_j(y'),
\end{align*}
which proves the second assertion. $\Box$

We note that $(g \ast P_{y'})(rx') = (g \ast P_{x'})(ry')$ and $\Lambda_t (g \ast P_{y'})(x) = (\Lambda_t g \ast
P_{y'})(x)$, these easy to prove formulae are often used in our proofs.

In this section $f_{m,y}$ stands for the harmonic function $f_{m, y}(x) = Q_m(x, y)$, $y \in \mathbb B$. We often
write $f_y$ instead of $f_{m, y}$. Let us collect some norm estimates of $f_y$.

\begin{lem}\label{est} For $0< p \leq \infty$ and $m>0$ we have
\begin{align}
M_\infty (f_{m, y}, r) & \leq C (1-|y|r)^{-n-m}, \qquad m \in \mathbb N. \label{fm0}\\
M_1(f_{m, y}, r) & \leq C (1-|y|r)^{-1-m}. \label{fm1}\\
\| f_{m, y} \|_{B^{p,1}_\alpha} & \leq C (1-|y|)^{\alpha - 1 - m}, \qquad m>\alpha - 1, \quad \alpha > 0. \label{fm2}\\
\| f_{m, y} \|_{B^{p,\infty}_\alpha} & \leq C (1-|y|)^{\alpha - n - m}, \qquad m \in \mathbb N, \quad
m > \alpha - n, \quad \alpha > 0. \label{fm2i}\\
\| f_{m, y} \|_{A^1_\alpha} & \leq C (1-|y|)^{\alpha - m}, \qquad m > \alpha > -1. \label{fm3} \\
\| f_{m, y} \|_{H^1_\alpha} & \leq C (1-|y|)^{\alpha - 1 - m} \qquad m > \alpha - 1, \quad \alpha \geq 0.\label{fm4}
\end{align}
\end{lem}

{\it Proof.} Using Lemma \ref{DSlemma} we obtain
\begin{equation*}
M_\infty(f_{m, y}, r) = \max_{x' \in \mathbb S} |Q_m(y, rx')| \leq \max_{x' \in \mathbb S}
\frac{C}{|\rho r x' - y'|^{n+m}} = C(1-r|y|)^{-n-m},
\end{equation*}
which gives (\ref{fm0}). The estimate (\ref{fm1}) follows from Lemma \ref{DSlemma}. The estimates (\ref{fm2}), for
finite $p$, and (\ref{fm3}) follow from Lemma \ref{rro} and (\ref{fm1}). Similarly, for finite $p$ (\ref{fm2i}) follows from (\ref{fm0}) and Lemma \ref{rro}. Next, using (\ref{fm1}),
$$\| f_{m, y} \|_{H^1_\alpha} \leq C \sup_{0 \leq r < 1} (1-r)^\alpha (1-r\rho)^{-m-1}, \qquad \rho = |y|.$$
The function $\phi(r) = (1-r)^\alpha (1-r\rho)^{-m-1}$ attains its maximum on $[0,1]$ at
$$r_0 = 1- (1-\rho) \frac{\alpha}{\rho(1+m-\alpha)},$$
as is readily seen by a simple calculus, and this suffices to establish (\ref{fm4}) and therefore (\ref{fm2}) for
$p = \infty$. Finally, (\ref{fm2i}) directly follows from Lemma \ref{DSlemma}. $\Box$

In this section we are looking for sufficient and/or necessary condition for a double indexed sequence $c$ to be
in $M_H(X, Y)$, for certain spaces $X$ and $Y$ of harmonic functions. We associate to such a sequence $c$ a
harmonic function
\begin{equation}\label{gc}
g_c(x) = g(x) = \sum_{k\geq 0} r^k \sum_{j=1}^{d_k} c_k^j Y^{(k)}_j(x'), \qquad x = rx' \in \mathbb B,
\end{equation}
and express our conditions in terms of $g_c$. Our main results give conditions in terms of fractional derivatives of
$g_c$, however it is possible to obtain some results on the basis of the following formula, contained in Lemma \ref{intcon}:
\begin{equation}\label{simple}
(c \ast f)(r^2x') = \int_{\mathbb S} (g \ast P_{y'})(rx') f(ry') dy'.
\end{equation}
Using continuous form of Minkowski's inequality, or more generally Young's inequality, this formula immediately gives the following proposition.

\begin{prop}
Let $c = \{ c_k^j : k \geq 0, 1 \leq j \leq d_k \}$ be a double indexed sequence and let $g(x) = \sum_{k\geq 0} r^k \sum_{j=1}^{d_k} c_k^j Y^{(k)}_j(x')$ be the corresponding harmonic function. If
$$\int_{\mathbb S} |(g \ast P_{y'})(rx')|^p dx' \leq C, \qquad y'\in \mathbb S, \quad 0 \leq r < 1,$$
then $c \in M_H(H^1, H^p)$.

More generally, if $1/q + 1/p = 1 + 1/r$, where $1 \leq p,q,r \leq \infty$, $\alpha + \gamma = \beta$, $\alpha,
\beta, \gamma \geq 0$ and $g \in H^p_\gamma$, then $c \in M_H(H^q_\alpha, H^r_\beta)$.
\end{prop}

The first part of the following lemma, which gives necessary conditions for $c$ to be a multiplier, is based on \cite{AS}.

\begin{lem}\label{nec}
Let $0 < p,q \leq \infty$, $1 \leq s \leq \infty$ and $m > \alpha - 1$. Assume a double indexed sequence
$c = \{ c_k^j : k \geq 0, 1 \leq j \leq d_k \}$ is a multiplier from $B^{p,1}_\alpha$ to $B^{q,s}_\beta$ and $g = g_c$
is defined in {\rm(\ref{gc})}. Then the following condition is satisfied:
\begin{equation}\label{mult}
N_s(g) = \sup_{0\leq\rho < 1} \sup_{y' \in \mathbb S} (1-\rho)^{m+1-\alpha + \beta} \left( \int_{\mathbb S}
|\Lambda_{m+1}(g \ast P_{x'})(\rho y')|^s dx'\right)^{1/s} < \infty,
\end{equation}
where the case $s = \infty$ requires usual modification.

Also, let $0 < p \leq \infty$, $1 \leq s \leq \infty$ and $m > \alpha - 1$. If a double indexed sequence
$c = \{ c_k^j : k \geq 0, 1 \leq j \leq d_k \}$ is a multiplier from $B^{p,1}_\alpha$ to $H^s_\beta$, then the above
function $g$ satisfies condition \rm{(\ref{mult})}.
\end{lem}

{\it Proof.}  Let $c \in M_H(B^{p,1}_\alpha, B^{q,s}_\beta)$, and assume both $p$ and $q$ are finite, the infinite cases require only small modifications. We have $\| M_c f\|_{B^{q,s}_\beta} \leq C \| f \|_{B^{p,1}_\alpha}$ for  $f$ in $B^{p,1}_\alpha$. Set $h_y = M_c f_y$, then we have
\begin{equation}\label{hay}
h_y(x) = \sum_{k\geq 0} r^k \rho^k \sum_{j=1}^k \frac{\Gamma(k+n/2+m+1)}{\Gamma(k+n/2)\Gamma(m+1)} c_k^j
Y^{(k)}_j(y') Y^{(k)}_j(x'), \quad x = rx' \in \mathbb B,
\end{equation}
moreover
\begin{equation}\label{cont}
\| h_y \|_{B^{q,s}_\beta} \leq C \| f_y \|_{B^{p,1}_\alpha}.
\end{equation}
This estimate and Lemma \ref{nec} give
\begin{equation}\label{hynorm}
\| h_y \|_{B^{q,s}_\beta} \leq C (1-|y|)^{\alpha - m - 1}, \qquad y \in \mathbb B.
\end{equation}
Note that $h_y(x) = \Lambda_{m+1}(g \ast P_{y'})(\rho x)$, using monotonicity of $M_s(h_y, r)$ we obtain:
\begin{align}\label{muka}
I_{y'}(\rho^2) & = \left( \int_{\mathbb S} |\Lambda_{m+1}(g \ast P_{x'})(\rho^2 y')|^s dx' \right)^{1/s}  =
\left( \int_{\rho}^1 (1-r)^{\beta q -1} r^{n-1} dr \right)^{-1/q} \notag\\
& \phantom{=} \times \left( \int_\rho^1 (1-r)^{\beta q -1} r^{n-1} \left( \int_{\mathbb S} |\Lambda_{m+1}
(g \ast P_{y'})(\rho^2 x')|^s dx'\right)^{q/s} dr\right)^{1/q} \notag \\
& \leq C(1-\rho)^{-\beta} \left( \int_\rho^1 (1-r)^{\beta q -1} r^{n-1} M_s^q(h_y, r) dr \right)^{1/q} \notag\\
& \leq C(1-\rho)^{-\beta} \| h_y \|_{B^{q,s}_\beta}.
\end{align}
Combining (\ref{muka}) and (\ref{hynorm}) we obtain
\begin{equation*}
\left( \int_{\mathbb S} |\Lambda_{m+1}(g \ast P_{x'})(\rho^2 y')|^s dx'\right)^{1/s} \leq C
(1 - \rho)^{\alpha - \beta - m -1},
\end{equation*}
which is equivalent to (\ref{mult}). The case $s = \infty$ is treated similarly.

Next we consider $c \in M_H(B^{p,1}_\alpha, H^s_\beta)$, assuming $0<p \leq \infty$.
Set $h_y = M_c h_y = g \ast f_y$. We have, by Lemma \ref{est},
$$\| f_y \|_{B^{p,1}_\alpha} \leq C(1-|y|)^{\alpha - m -1}, \qquad y \in \mathbb B,$$
and, by continuity of $M_c$, $\| h_y \|_{H^s_\beta} \leq C \| f_y \|_{B^{p,1}_\alpha}$. Therefore
$$\| h_y \|_{H^s_\beta} \leq C(1-|y|)^{\alpha - m -1}, \qquad y \in \mathbb B.$$
Setting $y = \rho y'$ we have
\begin{align*}
I_{y'}(\rho^2) & = \left( \int_{\mathbb S} |\Lambda_{m+1} (g \ast P_{x'})(\rho^2 y')|^s dx' \right)^{1/s}
= \left( \int_{\mathbb S} |\Lambda_{m+1} (g \ast P_y)(\rho x')|^s dx' \right)^{1/s}\\
& = M_s(h_y, \rho) \leq (1-|y|)^{-\beta} \| h_y \|_{H^s_\beta}.
\end{align*}
The last two estimates yield
$$\left( \int_{\mathbb S} |\Lambda_{m+1} (g \ast P_{x'})(\rho^2 y')|^s dx' \right)^{1/s} \leq C (1-|y|)^{\alpha - \beta
- m - 1},\qquad |y| = \rho$$
which is equivalent to (\ref{mult}). $\Box$

One of main results of this paper is a characterization of the space $M_H(B^{p,1}_\alpha, B^{q,1}_\beta)$ for $0 < p \leq q \leq \infty$. The following theorem treats the case $p>1$, while
Theorem \ref{main1} below covers the case $0 < p \leq 1$.

\begin{thm}\label{bpbp}
Let $1 < p \leq q \leq \infty$, and $m > \alpha - 1$. Then for a double indexed sequence
$c = \{ c_k^j : k \geq 0, 1 \leq j \leq d_k \}$ the following conditions are equivalent:

1. $c \in M_H(B^{p,1}_\alpha, B^{q,1}_\beta)$.

2. The function $g(x) = \sum_{k\geq 0} r^k \sum_{j=1}^{d_k} c_k^j Y^{(k)}_j(x')$ is harmonic in $\mathbb B$ and
satisfies the following condition
\begin{equation}\label{ng}
N_1(g)< \infty.
\end{equation}
\end{thm}

{\it Proof.} Since necessity of (\ref{ng}) is contained in Lemma \ref{nec} we prove sufficiency of condition (\ref{ng}). We assume $p$ and $q$ are finite, the remaining cases can be treated in a similar manner. Take $f \in B^{p, 1}_\alpha$ and set $h = M_c f$. Applying the operator $\Lambda_{m+1}$ to both sides of equation (\ref{simple}) we obtain
\begin{equation}\label{lamha}
\Lambda_{m+1} h(rx) = \int_{\mathbb S} \Lambda_{m+1} (g \ast P_{y'})(x) f(ry') dy'.
\end{equation}
Now we estimate the $L^1$ norm of the above function on $|x| = r$:
\begin{align}\label{form4}
M_1(\Lambda_{m+1} h, r^2) & \leq \int_{\mathbb S} M_1(\Lambda_{m+1} (g \ast P_{y'}), r) |f(ry')| dy' \notag \\
& \leq M_1(f, r) \sup_{y' \in \mathbb S} \int_{\mathbb S} | \Lambda_{m+1} (g \ast P_{y'})(rx')| dx' \notag \\
& \leq M_1(f, r)N_1(g) (1-r)^{\alpha - \beta -m -1}.
\end{align}
Since,
$$\int_0^1 M_1^p(h, r^2)(1-r)^{\beta p -1} r^{n-1} dr  \leq C \int_0^1 (1-r)^{p(m+1)} M_1^p(\Lambda_{m+1}h, r^2) (1-r)^{\beta p -1}r^{n-1} dr,$$
see \cite{DS}, we have
\begin{align*}
\| h \|_{B^{p,1}_\beta}^p &  \leq C \int_0^1 (1-r)^{p(m+1)} M_1^p(\Lambda_{m+1}h, r^2) (1-r)^{\beta p -1}r^{n-1} dr \\
& \leq C N_1^p(g) \int_0^1 M_1^p(f, r)(1-r)^{\alpha p - 1} r^{n-1} dr \\
& = CN_1^p(g) \| f \|_{B^{p,1}_\alpha}^p,
\end{align*}
and therefore $\| h \|_{B^{p,1}_\beta} \leq \| f \|_{B^{p,1}_\alpha}$. Since $\| h \|_{B^{q,1}_\beta} \leq C
\| h \|_{B^{p,1}_\beta}$ the proof is complete. $\Box$

Next we consider multipliers from $B^{p,1}_\alpha$ to $H^s_\beta$, in the case $0<p\leq 1$ we obtain a
characterization of the corresponding space.

\begin{thm}\label{bh}
Let $\beta \geq 0$, $0 < p \leq 1$, $s \geq 1$ and $m > \alpha - 1$. Then, for a double indexed sequence 
$c = \{ c_k^j : k \geq 0, 1 \leq j \leq d_k \}$ the following two conditions are equivalent:

1. $c \in M_H(B^{p,1}_\alpha, H^s_\beta)$.

2. The function $g(x) = \sum_{k\geq 0} r^k \sum_{j=1}^{d_k} c_k^j Y^{(k)}_j(x')$ is harmonic in $\mathbb B$ and
satisfies the following condition:
\begin{equation}\label{ng2}
N_s(g) < \infty.
\end{equation}
\end{thm}

{\it Proof.} The necessity of condition (\ref{ng2}) is contained in Lemma \ref{nec}. Now we turn to the sufficiency of (\ref{ng2}). We chose $f \in B^{p,1}_\alpha$ and set $h = c \ast f$. Then, by Lemma \ref{intcon}:
\begin{equation}\label{tool}
h(r^2 x') = 2 \int_0^1 \int_{\mathbb S} \Lambda_{m+1} (g \ast P_{\xi})(rR x') f(r R \xi) (1-R^2)^m R^{n-1} d\xi dR
\end{equation}
and this allows us to obtain the following estimate:
\begin{align*}
M_s(h, r^2) & \leq 2 \int_0^1 (1-R^2)^m R^{n-1} \left\| \int_{\mathbb S} \Lambda_{m+1} (g \ast P_\xi)(rRx') f(rR\xi)
d \xi \right\|_{L^s(\mathbb S, dx')} dR \notag\\
& \leq 2 \int_0^1 (1-R^2)^m R^{n-1} M_1(f, rR) \sup_{\xi \in \mathbb S} \|
\Lambda_{m+1} (g \ast P_\xi)(rRx') \|_{L^s} dR \notag \\
& \leq C N_s(g) \int_0^1 (1-R)^m M_1(f, rR) (1-rR)^{\alpha - \beta - m - 1} dR \notag \\
& \leq C N_s(g) \int_0^1 M_1(f, rR) (1-rR)^{\alpha - \beta - 1} dR \notag \\
& \leq C N_s(g) \int_0^1 M_1(f, rR) \frac{(1-R)^\alpha}{(1-rR)^{\beta + 1}} dR.
\end{align*}
Note that $M_1(f, rR)$ is increasing in $0\leq R < 1$, therefore we can combine Lemma \ref{wellkn} and the
above estimate to obtain, for $1/2 \leq r < 1$:
\begin{align*}
M_s^p(h, r^2) & \leq C N_s^p(g) \int_0^1 M_1^p(f, rR) \frac{(1-R)^{\alpha p + p -1}}{(1-rR)^{p\beta + p}} dR  \\
& \leq C N_s^p(g) (1-r)^{-p\beta} \int_0^1 M_1(f, R) (1-R)^{\alpha p -1} dR\\
& \leq C N_s^p(g) (1-r)^{-p\beta} \| f \|_{B^{p,1}_\alpha}^p \\
\end{align*}
Therefore $M_s(h, r^2) \leq C N_s(g) (1-r)^{-\beta} \| f \|_{B^{p,1}_\alpha}$, which completes the proof of the
Theorem. $\Box$

The omitted case $p = \infty$ is treated in our next theorem, which gives a characterization of the space $M_H(H^1_\alpha, H^p_\beta)$.
\begin{thm}\label{haha}
Let $\alpha \geq 0$, $\beta > 0$, $1 \leq p \leq \infty$ and $m > \alpha - 1$. Then for a double indexed sequence
$c = \{ c_k^j : k \geq 0, 1 \leq j \leq d_k \}$ the following conditions are equivalent:

1. $c \in M_H(H^1_\alpha, H^p_\beta)$.

2. The function $g(x) = \sum_{k\geq 0} r^k \sum_{j=1}^{d_k} c_k^j Y^{(k)}_j(x')$ is harmonic in $\mathbb B$ and
satisfies the following condition:
\begin{equation}\label{ng3}
N_p(g) < \infty.
\end{equation}
\end{thm}

In the case $p = \infty$ condition (\ref{ng3}) is interpreted in the usual manner.

{\it Proof.} Let us assume $c \in M_H(H^1_\alpha, H^p_\beta)$ and set $h_y = M_c f_y$ for $y \in \mathbb B$. Then we
have, by continuity of $M_c$ and by Lemma \ref{est}:
$$\| h_y \|_{H^p_\beta} \leq C \| f_y \|_{H^1_\alpha} \leq C (1-|y|)^{\alpha - m -1}.$$
On the other hand,
\begin{equation}
\| h_y \|_{H^p_\beta}  \geq (1-\rho)^\beta M_p(h_y, \rho) \geq (1-\rho)^\beta \left( \int_{\mathbb S} |\Lambda_{m+1}
(g \ast P_{x'})(\rho^2 y)|^p dx' \right)^{1/p},
\end{equation}
and the above estimates imply (\ref{ng3}). Now we prove sufficiency of the condition (\ref{ng3}). Choose
$f \in H^1_\alpha$ and set $h = c \ast f$. We apply continuous form of Minkowski's inequality to (\ref{lamha}) to
obtain
\begin{align*}
M_p(\Lambda_{m+1}h, r^2) & \leq M_1(f, r) \sup_{y' \in \mathbb S} M_p(\Lambda_{m+1} (g \ast P_{y'}), r)\\
& \leq N_p(g)(1-r)^{\alpha - \beta - m - 1} M_1(f, r).
\end{align*}
Therefore $\sup_{r<1} (1-r)^{m+1+\beta} M_p(\Lambda_{m+1}h, r) \leq C \| f \|_{H^1_\alpha}$ and it follows, see
\cite{DS}, that $\sup_{r<1} (1-r)^\beta M_p(h, r) \leq C \| f \|_{H^1_\alpha}$ as required. The case $p = \infty$
is treated the same way. $\Box$

Since $H^\infty_\beta = A^\infty_\beta$ the $p = \infty$ case of this theorem gives a complete description of the space $M_H(H^1_\alpha, A^\infty_\beta)$. The next proposition gives necessary conditions for $c$ to be in $M_H(X, A^\infty_\beta)$ for
some spaces $X$.

\begin{prop}
Let $m \in \mathbb N$ and $m > \alpha$. Let us consider, for a double indexed sequence $c = \{ c_k^j : k \geq 0, 1 \leq j \leq d_k \}$ the following conditions:

1. $c \in M_H(A^1_\alpha, A^\infty_\beta)$.

2. $c \in M_H(B^{p,1}_\alpha, A^\infty_\beta)$.

3. The function $g(x) = \sum_{k\geq 0} r^k \sum_{j=1}^{d_k} c_k^j Y^{(k)}_j(x')$ is harmonic in $\mathbb B$ and
satisfies the following condition:
\begin{equation}\label{multinhi}
M_t (g) = \sup_{0\leq\rho < 1} \sup_{x', y' \in \mathbb S} (1-\rho)^t
|\Lambda_{m+1}(g \ast P_{x'})(\rho y')| < \infty.
\end{equation}
Then we have: $1 \Rightarrow 3$ with $t = m + \beta - \alpha$ and
$2 \Rightarrow 3$ with $t = m+1+\beta-\alpha$.
\end{prop}

{\it Proof.} Let $X$ be one of the spaces $A^1_\alpha$, $B^{p,1}_\alpha$. As in the previous theorems, we choose a multiplier $c$ from $X$ to $A^\infty_\beta$ and note that $\| c \ast f \|_{A^\infty_\beta} \leq C \| f \|_X$.
We apply this to $f_y$, $y = \rho y' \in \mathbb B$ to obtain, with $h_y = c \ast f_y$, the estimate
$$\| h_y \|_{A^\infty_\beta} \leq C \| f_y \|_X.$$
Next,
\begin{align*}
\| h_y \|_{A^\infty_\beta} & \geq (1-\rho)^\beta M_\infty (h_y, \rho) = (1-\rho)^\beta \sup_{x' \in \mathbb S}
|h_y(\rho x')|\\
& = (1-\rho)^\beta \sup_{x' \in \mathbb S} |\Lambda_{m+1} (g \ast P_{x'})(\rho^2 y')|.\\
\end{align*}
and both implications follow from Lemma \ref{est}. $\Box$

The next theorem complements Theorem \ref{bpbp}, in less general form it appeared in \cite{AS} and for completness of
exposition we present, with permission of the authors, a proof.

\begin{thm}\label{main1}
Let $0<p\leq 1$, $m > \alpha - 1$ and $p \leq q \leq \infty$. Then for a double indexed sequence $c = \{ c_k^j : k \geq 0, 1 \leq j \leq d_k \}$ the following conditions are equivalent:

1. $c \in M_H(B^{p,1}_\alpha, B^{q,1}_\beta)$.

2. The function $g(x) = \sum_{k\geq 0} r^k \sum_{j=1}^{d_k} c_k^j Y^{(k)}_j(x')$ is harmonic in $\mathbb B$ and
satisfies the following condition
\begin{equation}\label{ng4}
N_1(g) < \infty.
\end{equation}
\end{thm}

{\it Proof.} Necessity of condition (\ref{ng4}) follows from Lemma \ref{nec}. Now we prove sufficiency of condition (\ref{ng4}). Let $f \in B^{p,1}_\alpha (\mathbb B)$ and set $h = c \ast f$. Then, using Lemma \ref{intcon}, we have:
\begin{align*}
\int_{\mathbb S} |h(r\rho x')| dx' & \leq \int_0^1 \int_{\mathbb S} \int_{\mathbb S} |\Lambda_{m+1} (g \ast P_{x'})
(rR\xi)| |f(\rho R \xi)|(1-R^2)^m R^{n-1} d\xi dx' dR \\
& \leq C \int_0^1 \left( \sup_{\xi \in \mathbb S} \int_{\mathbb S} |\Lambda_{m+1}(g \ast P_{x'})(rR\xi)| dx' \right)
\int_{\mathbb S} |f(\rho R\xi)| d \xi \\
& \phantom{\leq\leq} (1-R^2)^m R^{n-1} dR,
\end{align*}
and letting $\rho \to 1$ this gives
\begin{align*}
\int_{\mathbb S} |h(rx')| dx' &
\leq C \int_0^1 \left( \sup_{\xi \in \mathbb S} \int_{\mathbb S} |\Lambda_{m+1}(g \ast P_{x'})(rR\xi)| dx' \right)
\int_{\mathbb S} |f(R\xi)| d \xi \\
& \phantom{\leq\leq} (1-R^2)^m R^{n-1} dR.
\end{align*}
For each fixed $\xi \in \mathbb S$ the function
$$\psi_\xi(R) =  \int_{\mathbb S} |\Lambda_{m+1}(g \ast P_{x'})(rR\xi)| dx' = \int_{\mathbb S} |\Lambda_{m+1}
(g \ast P_\xi)(rRx')| dx'$$
is increasing for $0 \leq R < 1$, due to subharmonicity of $u_\xi(x) = |\Lambda_{m+1} (g \ast P_\xi)(rx)|$. Therefore
the function
$$G_r(R) = \left( \sup_{\xi \in \mathbb S} \int_{\mathbb S} |\Lambda_{m+1}(g \ast P_{x'})(rR\xi)| dx'\right)
\int_{\mathbb S} |f(R\xi)| d \xi, \qquad 0 \leq R < 1$$
is increasing and we can apply Lemma \ref{wellkn} to obtain
\begin{align*}
\left( \int_{\mathbb S} |h(rx'| dx'\right)^p & \leq C \left( \int_0^1 G_r(R) (1-R^2)^mR^{n-1}dR \right)^p \\
& \leq C \int_0^1 G_r(R)^p (1-R)^{mp + p -1} R^{n-1} dR.
\end{align*}
Since, for  $0 \leq r < 1$, $G_r(R) \leq N_1(g) M_1(f, R) (1-rR)^{\alpha - \beta -m -1}$ we have, using Lemma \ref{rro}
\begin{align*}
\| h \|_{B^{p,1}_\beta}^p & = \int_0^1 \left( \int_{\mathbb S} |h(rx')| dx' \right)^p
(1-r)^{p\beta - 1} r^{n-1} dr \\
& \leq C N_1(g)^p \int_0^1 M_1(f, R)^p (1-R)^{mp + p - 1} R^{n-1} \int_0^1
\frac{(1-r)^{p\beta -1} r^{n-1}dr}{(1-rR)^{p(m+1+\beta-\alpha)}} dR \\
& \leq C N_1(g)^p \int_0^1 M_1(f ,R)^p (1-R)^{p\alpha - 1} dR = C \| f \|_{B^{p,1}_\alpha}^p \\
\end{align*}
and we proved $\| h \|_{B^{p,1}_\beta} \leq C \| f \|_{B^{p,1}_\alpha}$ This, together with inequality
$\| h \|_{B^{q,1}_\beta} \leq C \| h \|_{B^{p,1}_\beta}$ finishes the proof. $\Box$

\section{Estimates for distances in harmonic function spaces in the unit ball and related problems in
$\mathbb R^{n+1}_+$}

In this section we investigate distance problems both in the case of the unit ball and in the case of the upper
half space.

\begin{lem}\label{emb}
Let $0 < p < \infty$ and $\alpha > -1$. Then there is a $C = C_{p, \alpha, n}$ such that for every $f \in A^p_\alpha
(\mathbb B)$ we have
$$|f(x)| \leq C (1-|x|)^{-\frac{\alpha + n}{p}} \| f \|_{A^p_\alpha}, \qquad x \in \mathbb B.$$
\end{lem}

{\it Proof.} We use subharmonic behavior of $|f|^p$ to obtain
\begin{align*}
|f(x)|^p & \leq \frac{C}{(1-|x|)^n} \int_{B(x, \frac{1-|x|}{2})} |f(y)|^p dy\\
& \leq C \frac{(1-|x|)^{-\alpha}}{(1-|x|)^n} \int_{B(x, \frac{1-|x|}{2})} |f(y)|^p (1-|y|)^\alpha dy\\
& \leq C(1-|x|)^{-\alpha - n} \| f \|_{A^p_\alpha}^p. \qquad \Box
\end{align*}
This lemma shows that $A^p_\alpha$ is continuously embedded in $A^\infty_{\frac{\alpha + n}{p}}$ and motivates the
distance problem that is investigated in Theorem \ref{thm4}.

\begin{lem}\label{lemb1}
Let $0<p<\infty$ and $\alpha > -1$. Then there is $C = C_{p, \alpha, n}$ such that for every $f \in \tilde A^p_\alpha$ and every $(x, t) \in \mathbb R^{n+1}_+$ we have
\begin{equation}\label{eqemb1}
|f(x, t)| \leq C y^{-\frac{\alpha + n + 1}{p}} \| f \|_{\tilde A^p_\alpha}.
\end{equation}
\end{lem}

The above lemma states that $\tilde A^p_\alpha$ is continuously embedded in $\tilde A^\infty_{\frac{\alpha + n +1}{p}}$, its proof is analogous to that of Lemma \ref{emb}.

For $\epsilon > 0$, $t > 0$ and $f \in h(\mathbb B)$ we set
$$U_{\epsilon, t}(f) = U_{\epsilon, t} =  \{ x \in \mathbb B : |f(x)|(1-|x|)^t \geq \epsilon \}.$$

\begin{thm}\label{thm4}
Let $p > 1$, $\alpha > -1$, $t = \frac{\alpha + n}{p}$ and $\beta > \max (\frac{\alpha + n}{p} - 1,
\frac{\alpha}{p})$, $\beta \in \mathbb N_0$. Set, for $f \in A^\infty_{\frac{\alpha + n}{p}}(\mathbb B)$:
$$t_1(f) = {\rm dist}_{A^\infty_{\frac{\alpha + n}{p}}} (f, A^p_\alpha),$$
$$t_2(f) = \inf \left\{ \epsilon > 0 : \int_{\mathbb B} \left( \int_{U_{\epsilon, t}}
|Q_\beta(x, y)|(1-|y|)^{\beta - t} dy \right)^p (1-|x|)^\alpha dx < \infty \right\}.$$
Then $t_1(f) \asymp t_2(f)$.
\end{thm}

{\it Proof.} We begin with inequality $t_1(f) \geq t_2(f)$. Assume $t_1(f) < t_2(f)$. Then there are
$0< \epsilon_1 < \epsilon$ and $f_1 \in A^p_\alpha$ such that $\|f-f_1\|_{A^\infty_t} \leq \epsilon_1$ and
$$\int_{\mathbb B} \left( \int_{U_{\epsilon, t}(f)} |Q_\beta (x, y)|(1-|y|)^{\beta - t}dy \right)^p
(1-|x|)^\alpha dx = +\infty.$$
Since $(1-|x|)^t |f_1(x)| \geq (1-|x|)^t|f(x)| - (1-|x|)^t|f(x) - f_1(x)|$ for every $x \in \mathbb B$ we conclude
that $(1-|x|)^t|f_1(x)| \geq (1-|x|)^t|f(x)| \geq (1-|x|)^t |f(x)| - \epsilon_1$ and therefore
$$(\epsilon - \epsilon_1) \chi_{U_{\epsilon, t}(f)}(x)(1-|x|)^{-t} \leq |f_1(x)|, \qquad x \in \mathbb B.$$
Hence
\begin{align*}
+\infty & = \int_{\mathbb B} \left( \int_{U_{\epsilon, t}(f)}|Q_\beta (x, y)|(1-|y|)^{\beta - t} dy \right)^p
(1-|x|)^\alpha dx\\
& = \int_{\mathbb B} \left( \int_{\mathbb B} \frac{\chi_{U_{\epsilon, t}(f)}(y)}{(1-|y|)^t}
|Q_\beta(x, y)|(1-|y|)^\beta dy \right)^p (1-|x|)^\alpha dx \\
& \leq C_{\epsilon, \epsilon_1} \int_{\mathbb B} \left( \int_{\mathbb B} |f_1(y)| |Q_\beta(x, y)|(1-|y|)^\beta dy \right)^p (1-|x|)^\alpha dx = M,
\end{align*}
and we are going to prove that $M$ is finite, arriving at a contradiction. Let $q$ be the exponent
conjugate to $p$. We have, using Lemma \ref{qbeta},
\begin{align*}
I(x) & = \left( \int_{\mathbb B} |f_1(y)| (1-|y|)^\beta |Q_\beta(x, y)| dy \right)^p \\
& = \left( \int_{\mathbb B} |f_1(y)| (1-|y|)^\beta |Q_\beta(x, y)|^{\frac{1}{n+\beta}(\frac{n}{p} + \beta -
\epsilon)} |Q_\beta(x, y)|^{\frac{1}{n+\beta}( \frac{n}{q} + \epsilon)} dy \right)^p \\
& \leq \int_{\mathbb B} |f_1(y)|^p (1-|y|)^{p\beta} |Q_\beta(x, y)|^{\frac{n + p\beta - p\epsilon}{n+\beta}} dy
\left( \int_{\mathbb B} |Q_\beta(x, y)|^{\frac{n+q\epsilon}{n+\beta}} dy \right)^{p/q}\\
& \leq C (1-|x|)^{-p\epsilon}
\int_{\mathbb B} |f_1(y)|^p (1-|y|)^{p\beta} |Q_\beta(x, y)|^{\frac{n + p\beta - p\epsilon}{n+\beta}} dy
\end{align*}
for every $\epsilon > 0$. Choosing $\epsilon > 0$ such that $\alpha - p \epsilon > -1$ we have, by Fubini's theorem
and Lemma \ref{qbeta}:
\begin{align*}
M & \leq C \int_{\mathbb B} |f_1(y)|^p (1-|y|)^{p\beta} \int_{\mathbb B} (1-|x|)^{\alpha - p\epsilon}
|Q_\beta(x, y)|^{\frac{n+p\beta - p\epsilon}{n+\beta}} dx dy \\
& \leq C \int_{\mathbb B} |f_1(y)|^p (1-|y|)^\alpha dy < \infty.
\end{align*}
In order to prove the remaining estimate $t_1(f) \leq C t_2(f)$ we fix $\epsilon > 0$ such that the integral appearing in the definition of $t_2(f)$ is finite and use Theorem \ref{intrep}, with $\beta > \max (t-1, 0)$:
\begin{align*}
f(x) & = \int_{\mathbb B \setminus U_{\epsilon, t}(f)} Q_\beta(x,y) f(y) (1-|y|^2)^\beta dy +
\int_{U_{\epsilon, t}(f)} Q_\beta(x, y) f(y) (1-|y|^2)^\beta dy \\
& = f_1(x) + f_2(x).
\end{align*}
Since, by Lemma \ref{qbeta}, $|f_1(x)|  \leq 2^\beta \int_{\mathbb B} |Q_\beta(x, y)| (1-|w|)^{\beta - t} dy  \leq C(1-|x|)^{-t}$ we have $\| f_1 \|_{A^\infty_t} \leq C \epsilon$. Thus it remains to show that $f_2 \in A^p_\alpha$
and this follows from
$$\| f_2 \|^p_{A^p_\alpha} \leq \| f \|_{A^\infty_t}^p \int_{\mathbb B} \left( \int_{U_{\epsilon, t}(f)}
|Q_\beta (x, y)| (1-|y|^2)^{\beta - t}  dy \right)^p (1-|x|)^\alpha dx < \infty. \quad \Box$$

The above theorem has a counterpart in the $\mathbb R^{n+1}_+$ setting. As a preparation for this result we need the
following analogue of Lemma \ref{qbeta}.

\begin{lem}\label{qm}
For $\delta > -1$, $\gamma > n + 1 + \delta$ and $m \in \mathbb N_0$ we have
$$\int_{\mathbb R^{n+1}_+} |Q_m(z, w)|^{\frac{\gamma}{n+m+1}} s^\delta dy ds \leq C t^{\delta - \gamma + n + 1},
\qquad t > 0.$$
\end{lem}

{\it Proof.} Using Fubini's theorem and estimate (\ref{estq}) we obtain
\begin{align*}
I(t) & = \int_{\mathbb R^{n+1}_+} |Q_m(z, w)|^{\frac{\gamma}{n+m+1}} s^\delta dy ds
\leq C \int_0^\infty s^\delta \left( \int_{\mathbb R^n} \frac{dy}{[|y|^2 + (s+t)^2]^\gamma} \right) ds \\
& = C\int_0^\infty s^\delta (s+t)^{n-\gamma} ds  = Ct^{\delta - \gamma + n + 1}. \qquad \Box
\end{align*}

For $\epsilon > 0$, $\lambda > 0$ and
$f \in h(\mathbb R^{n+1}_+)$ we set:
$$V_{\epsilon, \lambda}(f) = \{ (x, t) \in \mathbb R^{n+1}_+ : |f(x, t)|t^\lambda \geq \epsilon \}.$$

\begin{thm}
Let $p > 1$, $\alpha > -1$, $\lambda = \frac{\alpha + n + 1}{p}$, $m \in \mathbb N_0$ and
$m > \max (\frac{\alpha + n +1}{p} - 1, \frac{\alpha}{p})$. Set,
for $f \in \tilde A^\infty_{\frac{\alpha + n + 1}{p}}(\mathbb R^{n+1}_+)$:
$$s_1(f) = {\rm dist}_{\tilde A^\infty_{\frac{\alpha + n +1}{p}}} (f, \tilde A^p_\alpha),$$
$$s_2(f) = \inf \left\{ \epsilon > 0 : \int_{\mathbb R^{n+1}_+} \left( \int_{V_{\epsilon, \lambda}}
Q_m(z, w) s^{\beta - t}dy ds \right)^p t^\alpha dx dt < \infty \right\}.$$
Then $s_1(f) \asymp s_2(f)$.
\end{thm}

The proof of this theorem closely parallels the proof of the previous one, in fact, the role of Lemma \ref{qbeta}
is taken by Lemma \ref{qm} and the role of Theorem \ref{intrep} is taken by Theorem \ref{brthm}. We leave details
to the reader.

\end{document}